\documentclass[a4paper,11pt]{article}
\usepackage{amsthm}

\usepackage[T2A]{fontenc}
\usepackage[cp1251]{inputenc}
\usepackage[english,russian]{babel}
\usepackage[tbtags]{amsmath}
\usepackage{amsfonts,amssymb}

\usepackage{mathrsfs}
\usepackage{graphicx}

\usepackage{euscript,textcomp,verbatim,fancyhdr}
\usepackage[all]{xy}
\usepackage{epsfig}
\usepackage{latexsym}
\usepackage{mdwtab}           
\usepackage{graphics}

\theoremstyle{plain}
\newtheorem{Thm}{Теорема}

\newtheorem{Lem}{Лемма}

\theoremstyle{definition}
                             
                             \newtheorem{Def}{Определение}
                             \newtheorem*{Not}{Обозначение}

\newcommand{\lam}{\lambda}

\newcommand{\RR}{\mathbb R}

\newcommand{\NN}{\mathbb N}

\newcommand{\KK}{\mathbb{\widetilde K}}
\newcommand{\FF}{\mathbb{F}}
\newcommand{\DD}{\mathbb{D}}
\newcommand{\MM}{{\widetilde{\cal M}}}

\renewcommand{\P}{{\cal P}}

\newcommand{\D}{\mathscr D}

\newcommand{\N}{{\cal C}}

\newcommand{\eps}{\varepsilon}

\newcommand{\grad}{{\rm grad\,}}

\newcommand{\id}{{\rm id}}
\newcommand{\isot}{{\rm isot}}

\newcommand{\Diff}{{\rm Diff}}

\newcommand{\Ev}{{\rm Ev}}
\newcommand{\const}{{\mbox{\rm const}}}
\newcommand{\barc}{{\boldsymbol{c}}}
\newcommand{\J}{{\hat J}}
\newcommand{\conv}{{\mbox{\rm conv}}}
\renewcommand{\:}{{\colon\,}}

\renewcommand{\d}{\partial}
\newcommand{\Proof}{\noindent{\bf Доказательство. }}
\newcommand{\PProof}[1]{\noindent{\bf Доказательство {#1}. }}

\newcommand{\aapprox}{\sim} 
\renewcommand{\b}{y}
\renewcommand{\o}[1]{{\stackrel{{\scriptscriptstyle\circ}}{#1}}}
\renewcommand{\tilde}{\widetilde}
\renewcommand{\hat}{\widehat}
\renewcommand{\emptyset}{\varnothing}

\begin{document}
\date{}
\author{E.\,A.~Kudryavtseva}
\title{Special framed Morse functions on surfaces}

\maketitle

Let $M$ be a smooth closed orientable surface. Let $F$ be the space
of Morse functions on $M$,
and $\FF^1$ the space of framed Morse functions, both endowed with
$C^\infty$-topology. The space $\FF^0$ of special framed Morse
functions is defined. We prove that the inclusion mapping
$\FF^0\hookrightarrow\FF^1$ is a homotopy equivalence. In the case
when at least $\chi(M)+1$ critical points of each function of $F$ are
labeled, homotopy equivalences $\KK\sim\MM$ and
$F\sim\FF^0\sim\D^0\times\KK$ are proved, where $\KK$ is the complex
of framed Morse functions, $\MM\approx\FF^1/\D^0$ is the universal
moduli space of framed Morse functions, $\D^0$ is the group of
self-diffeomorphisms of $M$ homotopic to the identity.


{\bf MSC-class: } 58E05, 57M50, 58K65, 46M18

\bigskip
УДК 515.164.174+515.164.22+515.122.55

\begin{center}
{\LARGE Специальные оснащенные функции Морса на поверхностях} \bigskip\\
\large Е.\,А.~Кудрявцева
\bigskip
\end{center}

\begin{abstract}
Пусть $M$ --- гладкая замкнутая ориентируемая поверхность. Пусть $F$
-- пространство функций Морса на $M$,
и $\FF^1$ -- пространство оснащенных функций Морса, снабженные
$C^\infty$-топологией. Определено пространство $\FF^0$ специальных
оснащенных функций Морса и доказано, что отображение включения
$\FF^0\hookrightarrow\FF^1$ является
гомотопической эквивалентностью. В случае, когда у любой функции из
$F$ отмечено не менее чем $\chi(M)+1$ критических точек, доказаны
гомотопические эквивалентности $\KK\sim\MM$ и
$F\sim\FF^0\sim\D^0\times\KK$, где $\KK$ -- комплекс оснащенных
функций Морса, $\MM\approx\FF^1/\D^0$ -- универсальное пространство
модулей оснащенных функций Морса, $\D^0$ -- группа диффеоморфизмов
$M$, гомотопных тождественному.

\medskip
{\bf Ключевые слова: } функция Морса, оснащенная функция Морса,
комплекс оснащенных функций Морса, $C^\infty$-топология,
универсальное пространство модулей.
\end{abstract}

{\bf 1. Введение.} Настоящая работа является продолжением работ~\cite
{kp1,KP2} по изучению топологии пространства $F=F(M)$ функций
Морса на гладкой замкнутой ориентируемой поверхности $M$, снабженного
$C^\infty$-топологией.
В работе \cite {kp1} введено понятие оснащенной функции Морса (см.\
определение \ref {def:Morse}) и доказана гомотопическая
эквивалентность $F\sim\FF^1$ пространства $F$ функций Морса и
пространства $\FF^1=\FF^1(M)$ оснащенных функций Морса
(\cite{kp1,KP2}). В настоящей работе определено
пространство $\FF^0=\FF^0(M)$ специальных оснащенных функций Морса и
доказано, что отображение включения $\FF^0\hookrightarrow\FF^1$
является гомотопической эквивалентностью.
Согласно \cite {KP3,KP3addmz},
в большинстве случаев (см.\ (\ref {eq:main})) имеется гомотопическая
эквивалентность $\FF^1\sim R\times\MM$, где $R=R(M)$ -- одно из
многообразий $\RR P^3$,
$S^1\times S^1$ и точка (см.\ (\ref {eq:EE})), а $\MM=\MM(M)$ --
многообразие, гомеоморфное универсальному пространству модулей
оснащенных функций Морса. Мы доказываем, что
комплекс $\KK=\KK(M)$ оснащенных функций Морса (см.\ \cite {KP3,BaK})
является сильным деформационным ретрактом многообразия $\MM$.
Тем самым, мы доказываем гомотопическую эквивалентность $F\sim
R\times\KK$.

Обзор результатов о связных компонентах пространства функций Морса на
поверхности, коли\-чест\-ве и топологии его орбит при
действии диффеоморфизмов, о топологии про\-ст\-ранств функций с
умеренными особенностями,
связи с интегрируемыми системами имеется в~\cite {kp1} (см.\ также~\cite{F86}--\cite{MFSh88}).

\medskip
{\bf 2. Основные понятия и формулировка основного результата.}

 \begin{Def} \label {def:Morse}
Пусть $M$ --- гладкая (т.е.\ класса $C^\infty$) связная замкнутая
поверхность.

{\rm(A)} Пусть $F:=F_{p,q,r}(M)$ --- пространство функций Морса на
$M$, имеющих ровно $p$ критических точек локальных ми\-ни\-му\-мов,
$q$ седловых точек и $r$ точек локальных максимумов.
Для каждой функции $f\in F$ обозначим через $\N_{f,\lam}$ множество
ее критических точек индекса $\lam\in
\{0,1,2\}$, и $\N_f:=\N_{f,0}\cup\N_{f,1}\cup\N_{f,2}$.
Обозначим через $F^1$ подпространство в $F$, состоящее из функций
Морса $f\in F$, у которых все локальные минимумы равны $-1$, а все
локальные максимумы равны $1$.

(B) {\em Оснащенной функцией Морса на ориентированной поверхности
$M$} {\rm(см.~\cite[\S9]{kp1})} называется пара $(f,\alpha)$, где
$f\in F$ --- функция Морса на $M$, $\alpha$ --- замкнутая 1-форма на
$M\setminus(\N_{f,0}\cup\N_{f,2})$, такие что $2$-форма
$df\wedge\alpha$ не имеет нулей в $M\setminus\N_f$ и задает
положительную ориентацию, и в окрестности любой критической точки
$x\in\N_{f}$ существуют локальные координаты $u,v$, в которых либо
$f=u^2-v^2+f(x)$, $\alpha=d(2uv)$, либо
$f=\varkappa_{f,x}(u^2+v^2)+f(x)$,
$\alpha=\varkappa_{f,x}\frac{udv-vdu}{u^2+v^2}$, где
$\varkappa_{f,x}=\const\ne0$. Обозначим через
$\FF^1:=\FF^1_{p,q,r}(M)$ пространство оснащенных функций Морса
$(f,\alpha)$, таких
что $f\in F^1$. Снабдим $C^\infty$-топологией пространства $F,\FF^1$
(см.~\cite[\S4]{kp1}).
 \end{Def}

Из теоремы С.В.\ Матвеева (см.~\cite{Kmsb}) следует, что $\pi_0(F)=0$.

\begin{Not} [пермутоэдр и его грани]
{\it Пермутоэдр порядка $q\in\NN$} --- это выпуклый ($q-1$)-мерный
многогранник, вложенный в $q$-мерное пространство, вершины которого
получены пере\-ста\-нов\-ками координат вектора $(1,\dots,q)$.
Опишем его подробнее: пусть $e_1,\dots,e_q$ -- стандартный базис
$\RR^q$, и пусть $\P^{q-1}\subset\RR^q$ -- выпуклая оболочка
множества точек $P_\rho=\sum_{k=1}^q
\left(k-\frac{q+1}{2}\right)e_{\rho_k}$, $\rho\in\Sigma_q$.
Пермутоэдр $\P^{q-1}$
имеет ровно $q!$ вершин $P_\rho$, $\rho\in\Sigma_q$, а его
$(q-s)$-мерные грани находятся во взаимно однозначном соответствии с
упорядоченными разбиениями $J=(J_1,\dots,J_s)$ множества
$\{1,\dots,q\}$ на $s$ непустых подмножеств $J_1,\dots,J_s$.
А именно, грань $\tau_J$, отвечающая разбиению $J$ -- это выпуклая
оболочка множества точек
$(\Sigma_{r_1}\times\Sigma_{r_2-r_1}\times\ldots\times\Sigma_{r_s-r_{s-1}})(P_\rho)$,
где числа $0=r_0<r_1<\ldots<r_{s-1}<r_s=q$ и перестановка
$\rho\in\Sigma_q$ однозначно определяются условиями
 \begin{equation}\label{eq:J}
J_1=\{\rho_1,\dots,\rho_{r_1}\}, \
J_2=\{\rho_{r_1+1},\dots,\rho_{r_2}\}, \ \ldots, \
J_s=\{\rho_{r_{s-1}+1},\dots,\rho_{r_s}\},
 \end{equation}
$\rho_1<\ldots<\rho_{r_1}$, $\rho_{r_1+1}<\ldots<\rho_{r_2}, \ldots,
\rho_{r_{s-1}+1}<\ldots<\rho_{r_s}$. Здесь
$\Sigma_{r_1}\times\Sigma_{r_2-r_1}\times\ldots\times\Sigma_{r_s-r_{s-1}}$
-- подгруппа группы $\Sigma_q$, отвечающая разбиению
$\{1,\dots,q\}=\{1,\dots,r_1\}\sqcup\{r_1+1,\dots,r_2\}\sqcup\ldots\sqcup\{r_{s-1}+1,\dots,r_s\}$,
и действие перестановки $\sigma\in\Sigma_q$ на точке $P_\rho$ дает
точку $P_{\sigma\rho}$, где $(\sigma\rho)_i:=\rho_{\sigma_i}$, $1\le
i\le q$.
 Из описания граней многогранника $\P^{q-1}$ следует, что условие
$\tau_{\hat J}\subset\d\tau_J$ равносильно тому, что
разбиение $\hat J$ получается из разбиения $J$ путем измельчения.
\end{Not}

Для каждой функции $f\in F$ рассмотрим множество
$\N_{f,1}=:\{y_j\}_{j=1}^q\approx\{1,\dots,q\}$ ее седловых
критических точек (см.\ определение~\ref {def:Morse}(A)) и
евкли\-дово векторное пространство $0$-коцепей
 \begin{equation} \label{eq:H0f}
 H_f^0:=C^0(\N_{f,1};\RR)=\RR^{\N_{f,1}}\cong\RR^q.
 \end{equation}
Рассмотрим в пространстве $H^0_f$
многогранник $\P^{q-1}_f\subset H_f^0$, являющийся образом
пермутоэдра $\P^{q-1}\subset\RR^q$ при какой-либо биекции
$\N_{f,1}\to\{1,\dots,q\}$.
Рассмотрим ``вычисляющую'' 0-коцепь
 $$
 \barc=\barc(f):=f|_{\N_{f,1}}=(c_1,\dots,c_q)\in(-1;1)^{\N_{f,1}}\subset H_f^0,
 $$
сопоставляющую седловой точке $y_j\in \N_{f,1}$ значение
$c_j:=f(y_j)$, $1\le j\le q$. Сопоставим 0-коцепи
$\barc=(c_1,\dots,c_q)$
упорядоченное разбиение $J=J(\barc)=(J_1,\dots,J_s)$ множества седел
$\N_{f,1}\approx\{1,\dots,q\}$, определяемое свойствами~(\ref {eq:J})
и $c_{\rho_1}=\ldots=c_{\rho_{r_1}} <
 c_{\rho_{r_1+1}}=\ldots=c_{\rho_{r_2}} < \ldots <
c_{\rho_{r_{s-1}+1}}=\ldots=c_{\rho_{r_s}}$. (То есть, $J$ -- это
отношение частичного порядка на множестве седел $\N_{f,1}$ значениями
функции $f|_{\N_{f,1}}$.)

\begin{Def} [специальные оснащенные функции Морса] \label{def:sframed}

(A) {\it Сепаратрисой} оснащенной функции Морса $(f,\alpha)\in\FF$
назовем образ такой интегральной траектории $\gamma\colon\,(0,1)\to
M\setminus\N_f$ поля ядер 1-формы $\alpha$,
для которой оба предела $\lim_{t\to0^+}\gamma(t)$ и
$\lim_{t\to1^-}\gamma(t)$ существуют и принадлежат множеству $\N_f$,
и хотя бы один из этих пределов является седловой точкой.

(B) Оснащенную функцию Морса $(f,\alpha)\in\FF^1$ назовем {\it
специальной}, если либо $q=0$ (т.е. у функций из $F$ отсутствуют
седловые критические точки), либо выполнены следующие условия:
 \begin{enumerate}
 \item[(i)]
набор $\barc(f)\in H^0_f$ седловых критических значений принадлежит
многограннику $\frac2{q+1}\P^{q-1}_f$;
 \item[(ii)] пусть $\o\tau$ -- открытая грань многогранника
$\P^{q-1}_f$, содержащая точку $\frac{q+1}2\barc(f)$, и пусть
$J=(J_1,\dots,J_s)$ -- соответствующее упорядоченное разбиение
множества седел $\N_{f,1}$ на непустые подмножества, т.е.\
$\tau=\tau_J^{q-s}$; тогда для любых двух седел $\b_i,\b_j\in\N_{f,1}$ из
одного и того же подмножества $J_k$ разбиения ($1\le k\le s$)
не существует сепаратрисы, соединяющей $\b_i$ и $\b_j$.
 \end{enumerate}
Пусть $\FF^0:=\FF^0_{p,q,r}(M)$ -- пространство специальных
оснащенных функций Морса.
 \end{Def}

Группа диффеоморфизмов $\D^\pm:=\Diff(M)$ действует справа на пространстве $\FF^1$
очевидным образом {\rm(см.~\cite[обозначение 2.3]{kp1})}. Очевидно,
что $\FF^0$ является $\D^\pm$-инвариантным.

\begin{Thm}  \label{thm:FF0}
Пусть $M$ --- замкнутая ориентированная поверхность,
$F=F_{p,q,r}(M)$ -- пространство функций Морса на $M$,
$\FF^0\subset\FF^1$ --- соответствующие пространства оснащенных
функ\-ций Морса {\rm(см.\ опре\-де\-ле\-ния~$\ref {def:Morse}$, $\ref
{def:sframed}$)}. Отображение включения $\FF^0\hookrightarrow\FF^1$
является гомотопической экви\-ва\-лент\-ностью, причем
соответствующие отображения и гомотопии могут быть выбраны
$\D^\pm$-эквивариантными и сохраняющими отображение
$\FF^1\to M^{p+q+r}/\Sigma_{p+q+r}$,
$(f,\alpha)\mapsto\N_f$.
 \end{Thm}

{\bf 3. Доказательство теоремы \ref {thm:FF0}.}
%
Предположим, что число седел $q\ge1$. Фиксируем вещественное число $\kappa>0$.
Так как $\P^{q-1}\subset[-\frac{q-1}2;\frac{q-1}2]^q$, то
$\kappa\P^{q-1}\subset(-\frac{q\kappa}2;\frac{q\kappa}2)^q$.
Рассмотрим отображение
 $$
\pi_\kappa\colon\,\RR^q\to\kappa\P^{q-1},
 $$
переводящее любую точку $\barc\in\RR^q$ в такую точку
$\barc'\in\kappa\P^{q-1}$, что $|\barc-\barc'|\le|\barc-\barc''|$ для
любой точки $\barc''\in\kappa\P^{q-1}$. В силу выпуклости пермутоэдра $\P^{q-1}$,
такое отображение единственно.

Фиксируем оснащенную функцию Морса $(f,\alpha)\in\FF^1$. Выберем
какую-нибудь нумерацию $\N_{f,1}\approx\{1,\dots,q\}$ множества
седловых точек, и с ее помощью отождествим $H^0_f\cong\RR^q$ и
$\P^{q-1}_f\cong\P^{q-1}$. Пусть
$\pi_{\kappa,f}\:H^0_f\to\kappa\P^{q-1}_f$ -- композиция
$H^0_f\cong\RR^q\stackrel{\pi_\kappa}{\longrightarrow}\kappa\P^{q-1}\cong\kappa\P^{q-1}_f$.
Рассмотрим на $M\setminus(\N_{f,0}\cup\N_{f,2})$ гладкое поле
неотрицательно определенных в каждой точке квадратичных форм $(df)^2+\alpha^2$.
Аналогично римановым метрикам, этому полю квадратичных форм отвечает
функция длины $L(\gamma)$ регулярных кусочно-гладких путей
$\gamma$ на $M\setminus(\N_{f,0}\cup\N_{f,2})$.
Определим расстояние $\rho(x,y)=\rho_{f,\alpha}(x,y):=\inf(L(\gamma))$, $x,y\in M\setminus(\N_{f,0}\cup\N_{f,2})$, где нижняя
грань берется по регулярным кусочно-гладким путям $\gamma$ на
поверхности $M\setminus(\N_{f,0}\cup\N_{f,2})$ из $x$ в $y$.
Определим вещественное число $d_{j_1,j_2}=d_{j_2,j_1}$ равным
расстоянию $\rho_{f,\alpha}(\b_{j_1},\b_{j_2})$
между седловыми точками $\b_{j_1}$ и $\b_{j_2}$ при $1\le j_1<j_2\le
q$. Пусть $\barc(f)=(c_1,\dots,c_q)\in(-1;1)^q$ -- набор седловых значений.
Положим
 \begin{equation}\label{eq:eps:f:alpha}
 \eps=\eps(f,\alpha) :=\frac13\min\left\{1,\ \min_{1\le j_1<j_2\le q}d_{j_1,j_2},\
 1-\max_{1\le j\le q}|c_j|
\right\} \in \left(0;\frac13\right],
 \end{equation}
 \begin{equation} \label {eq:c'}
 \barc'=(c'_1,\dots,c'_q):=\pi_{\frac{\eps}{q+1},f}(\barc(f))
 \in\frac{\eps}{q+1}\P^{q-1}_f \subset\left(-\frac{\eps}2;\frac{\eps}2\right)^{\N_{f,1}}\subset (-1;1)^{\N_{f,1}}\cong(-1;1)^q.
 \end{equation}
Пусть $\o\tau$ -- открытая грань многогранника $\P^{q-1}_f$,
содержащая точку $\frac{q+1}{\eps}\barc'\in\P^{q-1}_f$. Согласно
обозна\-че\-нию,
эта грань имеет вид $\tau=\tau_\J^{q-s}$ для некоторого $s\in[1;q]$ и
некоторого разбиения $\J=(\J_1,\dots,\J_s)$
множества $\N_{f,1}\approx\{1,\dots,q\}$ на $s$ подмножеств
$\J_k\approx\{\rho_{r_{k-1}+1},\dots,\rho_{r_k}\}$, $1\le k\le s$,
см.~(\ref{eq:J}). Напомним, что
грань $\tau$ состоит из всех точек $\barc=(c_1,\dots,c_q)\in
H^0_f\cong\RR^q$, для которых
 \begin{equation} \label{eq:tau'}
(c_{\rho_{r_{k-1}+1}},\dots,c_{\rho_{r_k}})
 \in\conv\left(\Sigma_{r_k-r_{k-1}}\ (r_{k-1}+1-\frac{q+1}2,\dots,r_k-\frac{q+1}2)\right),
 \quad 1\le k\le s,
 \end{equation}
см.~(\ref{eq:J}).
Положим
 $t_0:=-1+\frac{\eps}2$,
 $t_k:=c_{\rho_{r_k}}(f)-c'_{\rho_{r_k}}$ при $1\le k\le s$,
 $t_{s+1}:=1-\frac{\eps}2$.
Определим функцию
$h=h_{\barc(f),\eps}\:[-\frac{\eps}2;\frac{\eps}2]\to[-1;1]$ формулой
 \begin{equation}\label {eq:h:c:eps}
 h(t)=h_{\barc(f),\eps}(t):= t + t_0
 + \sum_{k=0}^{s} (t_{k+1}-t_k) I_{(\frac{r_k}{q+1}-\frac12)\eps,(\frac{r_k+1}{q+1}-\frac12)\eps} (t),
 \quad -\frac{\eps}2\le t\le\frac{\eps}2,
 \end{equation}
где $I_{a,b}\in C^\infty(\RR)$ -- гладкое двупараметрическое
семейство функций с параметрами $a<b$, такое что $I_{a,b}'\ge0$,
$I_{a,b}|_{(-\infty;(2a+b)/3]}=0$ и $I_{a,b}|_{[(a+2b)/3;+\infty)}=1$
(определенное, например, как в \cite[(6)]{kp1}).

\begin{Lem}\label{lem:h}
Функция $h=h_{\barc(f),\eps}\:[-\frac{\eps}2;\frac{\eps}2]\to[-1;1]$
является диффеоморфизмом отрезков, причем $t_0<t_1\le\ldots\le
t_s<t_{s+1}$. Кроме того, выполнены следующие условия:
\begin{enumerate}
 \item[1)] $h_{\barc(f),\eps}(c_j')=c_j'+t_k=c_j(f)$ для любого $j\in
\J_k\subset\{1,\dots,q\}$, $1\le k\le s$;
 \item[2)] в некоторой окрестности множества
$\{c_1',\dots,c_q'\}\cup\{-\frac{\eps}2,\frac{\eps}2\}$ в
$[-\frac{\eps}2;\frac{\eps}2]$ выполнено $h'\equiv1$;
 \item[3)] если $1\le u\le s$, $0=k_0<k_1<\ldots<k_u=s$ и
 \begin{equation} \label{eq:perturb}
 \barc(f)\in\overline{\pi_{\frac{\eps}{q+1},f}^{-1}(\frac{\eps}{q+1}\o\tau^{q-u}_{
 (\J_1\cup\ldots\cup\J_{k_1},\J_{k_1+1}\cup\ldots\cup\J_{k_2},\ldots,\J_{k_{u-1}+1}\cup\ldots\cup\J_{k_u})
 })},
 \end{equation}
то $t_1=\ldots=t_{k_1}\le t_{k_1+1}=\ldots=t_{k_2}\le\ldots\le
t_{k_{u-1}+1}=\ldots=t_{k_u}$.
\end{enumerate}
\end{Lem}

\Proof Покажем, что $t_0<t_1$ и $t_s<t_{s+1}$. Так как
$c_j(f)\in[-1+3\eps;1-3\eps]$ в силу
(\ref{eq:eps:f:alpha}), $c_j'\in(-\frac{\eps}2;\frac{\eps}2)$ в
силу~(\ref {eq:c'}), то
$t_1-t_0=c_{\rho_{r_1}}(f)-c'_{\rho_{r_1}}+1-\frac{\eps}2>0$,
$t_{s+1}-t_s=1-\frac{\eps}2-c_{\rho_{r_s}}(f)+c'_{\rho_{r_s}}>0$.

Покажем, что $t_1\le\ldots\le t_s$. Из вида граней пермутоэдра
$\P^{q-1}$ (см.\ обозначение) следует, что
$$
 \frac{\eps}{q+1}\P^{q-1}_f=\left\{\left.\phantom{I^{I^I}_{j_j}}\!\!\!\!\!\!\!
 \barc=(c_1,\dots,c_q)\in H^0_f\cong\RR^q\ \right|
 \Phi(\barc)=0,\
 \Phi_{J_1}(\barc)\le 0,\ \emptyset\ne
 J_1\subsetneq\{1,\dots,q\}\right\},
$$
где линейные функции $\Phi,\Phi_{J_1}\:H^0_f\to\RR$ определены
формулами
 $$
\Phi(\barc):=c_1+\ldots+c_q, \quad
 \Phi_{J_1}(\barc)
 :=-\sum_{j\in J_1}c_j+\frac{\eps}{q+1}\sum_{i=1}^{|J_1|}(i-\frac{q+1}{2})
 =-\sum_{j\in J_1}c_j+\frac{|J_1|-q}{2(q+1)}\eps |J_1|.
 $$
В точке $\barc'$ достигается минимум функции
$\hat\Phi(\barc):=\frac12|\barc(f)-\barc|^2$ по всем точкам
$\barc\in\frac{\eps}{q+1}\P^{q-1}_f$. Равенство $\Phi_{J_1}(\barc)=0$
равносильно тому, что $\barc\in\tau_{(J_1,\overline {J_1})}^{q-2}$.
Поэтому равенство $\Phi_{J_1}(\barc')=0$ равносильно тому, что
$J_1=\hat J_1\cup\ldots\cup \hat J_k$ для некоторого $k\in[1;s-1]$.
По теореме Куна-Таккера~\cite{MT} выполнено
$$
 \grad\hat\Phi(\barc')+\lam\,\grad\Phi(\barc')
 +\sum_{k=1}^{s-1}
 \lam_k\,\grad\Phi_{\hat J_1\cup\ldots\cup \hat J_k}(\barc')=0
$$
для некоторых множителей Лагранжа
$\lam,\lam_1,\dots,\lam_{s-1}\in\RR$, $\lam_k\ge 0$. Отсюда
$$
 \barc(f)-\barc'=\lam(e_1+\ldots+e_q)-\sum_{k=1}^{s-1}\lam_k(e_{\rho_1}+\ldots+e_{\rho_{r_k}})
 =\sum_{k=1}^s(\lam-\sum_{i=k}^{s-1}\lam_i)(e_{\rho_{r_{k-1}+1}}+\ldots+e_{\rho_{r_k}}).
$$
Поэтому для любого $j\in\hat J_k$ имеем
$c_j(f)-c'_j=\lam-\sum_{i=k}^{s-1}\lam_i$.
Отсюда $t_1\le\ldots\le t_s$.

Из равенств $h(-\frac{\eps}2)=-\frac{\eps}2+t_0=-1$,
$h(\frac{\eps}2)=\frac{\eps}2+t_{s+1}=1$, неравенств
$t_0<t_1\le\ldots\le t_s<t_{s+1}$ и неубывания функций $I_{a,b}$
следует, что $h\in\Diff([-\frac{\eps}2;\frac{\eps}2],[-1;1])$.

Осталось доказать выполнение свойств 1)--3). Для каждой функции
$I_{a,b}$ из определения функции $h$ имеем $b-a=\frac{\eps}{q+1}$.
Отсюда и из определения функции $I_{a,b}$ следует, что
в $\frac{\eps}{3(q+1)}$\,--окрестности точки $-\frac{\eps}2$ в
$[-\frac{\eps}2;\frac{\eps}2]$ имеем $h(t)\equiv t+t_0$, в
$\frac{\eps}{3(q+1)}$\,--окрестности отрезка
$\left[(\frac{r_{k-1}+1}{q+1}-\frac12)\eps;(\frac{r_k}{q+1}-\frac12)\eps\right]$
имеем $h(t)\equiv t+t_k$, $1\le k\le s$, а в
$\frac{\eps}{3(q+1)}$\,--окрестности точки $\frac{\eps}2$ в
$[-\frac{\eps}2;\frac{\eps}2]$ имеем $h(t)\equiv t+t_{s+1}$.
С другой стороны, из условия $\frac{q+1}{\eps}\barc'\in\o\tau_{\hat
J}$ следует, что для любого $j\in\hat J_k$ выполнено
$c_j'\in\left[(\frac{r_{k-1}+1}{q+1}-\frac12)\eps;(\frac{r_k}{q+1}-\frac12)\eps\right]$,
см.~(\ref{eq:tau'}), откуда
$h(c'_j)=c'_j+t_k$, $1\le k\le s$.

Для любого
 $\barc\in\pi_{\frac{\eps}{q+1},f}^{-1}(\frac{\eps}{q+1}\o\tau_{
 (\J_1\cup\ldots\cup\J_{k_1},\J_{k_1+1}\cup\ldots\cup\J_{k_2},\ldots,\J_{k_{u-1}+1}\cup\ldots\cup\J_{k_u})})$
имеем $\barc-\pi_{\frac{\eps}{q+1},f}(\barc)
=\sum_{i=1}^u(\mu-\mu_i-\ldots-\mu_{u-1})(e_{\rho_{r_{k_{i-1}-1}+1}}+\ldots+e_{\rho_{r_{k_i}}})$
для некоторых $\mu,\mu_1,\dots,\mu_u\in\RR$, $\mu_i\ge 0$. Значит,
в случае~(\ref{eq:perturb}) выполнено
$t_{k_{i-1}+1}=\ldots=t_{k_i}=\mu-\mu_i-\ldots-\mu_{u-1}$, $1\le i\le
u$, где $k_0:=0$. Лемма \ref {lem:h} доказана. \qed

\medskip
\PProof{теоремы 1} {\it Шаг 1.} Если $q=0$, то $\FF^0=\FF^1$ и
доказывать нечего. Пусть далее число седел $q>0$. Сопоставим
оснащенной функции Морса $(f,\alpha)\in\FF^1$ функцию Морса $\tilde
f:=h^{-1}\circ f$, где $h=h_{\barc(f),\eps(f,\alpha)}$ --
диффеоморфизм из леммы~\ref {lem:h}.
Из леммы~\ref {lem:h} и определения~\ref {def:Morse} следует, что
$(\tilde f,\alpha)\in\FF$.
Так как $h[-1;1]=[-\frac{\eps}2;\frac{\eps}2]$, то
$\frac2{\eps}(\tilde f,\alpha) \in\FF^1$. Покажем, что $\frac2{\eps}(\tilde f,\alpha)\in\FF^0$.
Так как $\barc(\tilde f)=\barc'\in \frac{\eps}{q+1}\P^{q-1}_f$, то
$\barc(\frac2{\eps}\tilde f)=\frac2{\eps}\barc'\in
\frac2{q+1}\P^{q-1}_f$. Это доказывает выполнение условия~(i) из
определения \ref {def:sframed} пространства $\FF^0$. Осталось
проверить условие~(ii).
Пусть $\gamma$ -- сепаратриса оснащенной функции $(f,\alpha)$,
соединяющая седловые точки $\b_i,\b_j\in\N_{\tilde f,1}=\N_{f,1}$,
$1\le i<j\le q$. Тогда
 $$
 |c_i-c_j|
 = |c_i-c_j| + |\int_\gamma \alpha|
 = |\int_\gamma df|+|\int_\gamma \alpha|
 = \int_\gamma(|df|+|\alpha|)
 \ge \int_\gamma\sqrt{df^2+\alpha^2}
 \ge \rho_{f,\alpha}(\b_i,\b_j)
 \ge 3\eps
 $$
в силу (\ref{eq:eps:f:alpha}). Предположим также, что
седловые точки принадлежат одному и тому же подмножеству $\J_k$
разбиения $\J$, где $1\le k\le s$. Тогда из
леммы~\ref {lem:h} следует, что $c_i=c_i'+t_k$ и $c_j=c_j'+t_k$, а
потому $c_i-c_j=c_i'-c_j'$. Отсюда, с учетом
$c_i',c_j'\in(-\frac{\eps}2;\frac{\eps}2)$, получаем
$|c_i-c_j|=|c_i'-c_j'|<\eps$,
что противоречит вышеприведенному неравенству.
Таким образом, $\frac2{\eps}(\tilde f,\alpha)\in\FF^0$, и возникает
отображение
 \begin{equation}\label{eq:p3}
p_3\:\FF^1\to\FF^0, \quad (f,\alpha)\mapsto\frac2{\eps}(\tilde
f,\alpha)
 = \frac2{\eps(f,\alpha)}\left(h_{\barc(f),\eps(f,\alpha)}^{-1}\circ f,\alpha\right).
 \end{equation}

{\it Шаг 2.} Докажем непрерывность отображения $p_3\:\FF^1\to\FF^0$.

Для любых $\kappa\in(0;1)$ и
$\barc\in (-1+\kappa;1-\kappa)^q$ рассмотрим диффеоморфизмы
$m_{\frac\kappa2}\:[-1;1]\to[-\frac\kappa2;\frac\kappa2]$,
$t\mapsto\frac\kappa2t$, и $h_{\barc,\kappa}\circ m_{\frac\kappa2}\:[-1;1]\to[-1;1]$, см.~(\ref
{eq:h:c:eps}). Покажем, что сопоставление
 $$
 H\:
(\barc,\kappa) \mapsto h_{\barc,\kappa}\circ m_{\frac\kappa2}\in\Diff[-1;1]
 $$
непрерывно на множестве таких пар $(\barc,\kappa)$, что
$\barc\in (-1+\kappa;1-\kappa)^q$ и $\kappa\in(0;1)$. Пусть $(\barc,\kappa)$ --
любая пара из этого множества, и $\o\tau:=\o\tau_{(\J_1,\dots,\J_s)}$
-- открытая грань многогранника $\P^{q-1}$, содержащая точку
$\frac{q+1}{\kappa}\pi_{\frac{\kappa}{q+1}}(\barc)\in\P^{q-1}$. Из
определения $h_{\barc,\kappa}$ и непрерывности
$\pi_{\frac{\kappa}{q+1}}$ следует непрерывность ограничения $H$
на множество таких пар $(\barc',\kappa')$, что
$\frac{q+1}{\kappa'}\pi_{\frac{\kappa'}{q+1}}(\barc')\in\o\tau$.
Осталось проверить непрерывность в точке $(\barc,\kappa)$ ограничения
$H$ на множество таких пар $(\barc',\kappa')$, что
   $\barc'\in\pi_{\frac{\kappa'}{q+1}}^{-1}(\frac{\kappa'}{q+1}\o\tau)
 \cup\overline{\pi_{\frac{\kappa'}{q+1}}^{-1}(\frac{\kappa'}{q+1}\o\tau_1)}$,
для любой открытой грани $\o\tau_1$ многогранника $\P^{q-1}$. Из
непрерывности $\pi_{\frac{\kappa}{q+1}}$ следует, что если
$\barc\in\d\left(\pi_{\frac{\kappa}{q+1}}^{-1}(\frac{\kappa}{q+1}\o\tau_1)\right)$,
то $\tau\subset\tau_1$. Тогда грань $\tau_1$ имеет вид
$\tau_1=\tau_{(\J_1\cup\ldots\cup\J_{k_1},\J_{k_1+1}\cup\ldots\cup\J_{k_2},
\ldots,\J_{k_{u-1}+1}\cup\ldots\cup\J_{k_u})}$. Согласно
лемме~\ref{lem:h}, имеем $t_1=\ldots=t_{k_1}\le
t_{k_1+1}=\ldots=t_{k_2}\le\ldots\le t_{k_{u-1}+1}=\ldots=t_{k_u}$.
Значит, $h_{\barc,\kappa}(t) = t + t_0 + \sum_{i=0}^{u}
(t_{k_i+1}-t_{k_i})
I_{(\frac{r_{k_i}}{q+1}-\frac12)\kappa,(\frac{r_{k_i}+1}{q+1}-\frac12)\kappa}
(t)$, где $k_0:=0$, т.е.\ диффеоморфизм $h_{\barc,\kappa}$
определяется той же формулой, что и $h_{\barc',\kappa'}$ при
$\barc'\in\pi_{\frac{\kappa'}{q+1}}^{-1}(\frac{\kappa'}{q+1}\o\tau_1)$.
Непрерывность $H$ доказана.

Согласно~(\ref {eq:p3}), выполнено
 $p_3(f,\alpha)=\frac2{\eps(f,\alpha)}(h_{\barc(f),\eps(f,\alpha)}^{-1}\circ f,\alpha)
 = ((H(\barc(f),\eps(f,\alpha)))^{-1}\circ f,\frac2{\eps(f,\alpha)}\alpha)$.
Отображение $\Diff[-1;1]\to\Diff[-1;1]$, $H\mapsto H^{-1}$,
непрерывно в $C^\infty$-топологии в силу~\cite[лемма 10.1]{kp1}.
Отсюда, а также из непрерывности функции $\eps=\eps(f,\alpha)$,
сопоставления $f\mapsto\Sigma_q\barc(f)\in\RR^q/\Sigma_q$
и отображения $H$ (см.\ выше), следует непрерывность отображения
$p_3$.

{\it Шаг 3.}
Покажем, что отображение включения $i_3\:\FF^0\hookrightarrow\FF^1$
является гомотопической экви\-ва\-лент\-ностью.
Композиция $i_3\circ p_3\:\FF^1\to\FF^1$ переводит
$(f,\alpha)\mapsto\frac2{\eps}(h^{-1}\circ f,\alpha)$, где
$\eps=\eps(f,\alpha)$,
$h:=h_{\barc(f),\eps}\in\Diff([-\frac{\eps}2;\frac{\eps}2],[-1;1])$,
см.\ лемму~\ref {lem:h}.
Зададим гомотопию этой композиции в $\id_{\FF^1}$ формулой
 $$
(f,\alpha)\mapsto(f_t,\alpha_t):=(1-t)\cdot \frac2{\eps}(h^{-1}\circ
f,\alpha)+t\cdot(f,\alpha),
 \quad 0\le t\le1.
 $$
Функция $f_t$ является функцией Морса из $F^1$ с теми же линиями
уровня и тем же множеством $\N_f$ критических точек, что и у функции
$f$, причем в силу части~2 леммы~\ref{lem:h}
функция $f_t$ в окрестности каждой критической точки отличается от
$f$ домножением на положительное число $(1-t)\frac2{\eps}+t$ и
прибавлением некоторой константы,
а форма $\alpha_t$ отличается от $\alpha$ домножением на то же самое
число. Отсюда следует, что гомотопия не выводит из пространства
$\FF^1$.

Покажем, что ограничение на $\FF^0$ этой гомотопии не выводит из
пространства $\FF^0$.
При отображении
$m_{\frac2\eps}\circ\pi_{\frac\eps{q+1},f}|_{\frac{2}{q+1}\P^{q-1}_f}$
образ любой открытой грани $\frac{2}{q+1}\o\tau$ многогранника
$\frac{2}{q+1}\P^{q-1}_f$ лежит в грани $\frac{2}{q+1}\tau$. Отсюда
следует, что если $\barc(f)\in\frac{2}{q+1}\o\tau$, то
$\barc(f_t)=(1-t)\frac2\eps\pi_{\frac\eps{q+1},f}(\barc(f))+t\barc(f)\in\frac{2}{q+1}\o\tau$
при $0<t\le1$. Значит, паре $(f_t,\alpha_t)$ отвечает
то же разбиение $J=(J_1,\dots,J_s)$, что и паре $(f,\alpha)$. Поэтому
выполнение условия~(ii) для пар $(f_t,\alpha_t)$, $0<t\le1$, следует
из выполнения аналогичного условия для пары $(f,\alpha)$, с учетом
$\alpha_t=((1-t)\frac2{\eps}+t)\alpha$.

По построению, оба отображения $i_3,p_3$ и гомотопия
$\D^\pm$-эквивариантны. Теорема \ref {thm:FF0} доказана.

\medskip
{\bf 4. Применение к исследованию гомотопического типа пространства
функций Морса.}
Обозначим через $F:=F_{p,q,r;\hat p,\hat q,\hat r}(M)$
про\-стран\-ство, полученное из $F_{p,q,r}(M)$ введением нумерации у
некоторых из критических точек (называемых отмеченными) для функций
$f\in F_{p,q,r}(M)$, где $\hat p,\hat r,\hat q$ -- количества
отмеченных критических точек локальных минимумов, максимумов и
седловых точек соответственно. Пусть $\FF^0:=\FF^0_{p,q,r;\hat p,\hat
q,\hat r}(M)$ и $\FF^1:=\FF^1_{p,q,r;\hat p,\hat q,\hat r}(M)$ --
соответствующие пространства оснащенных функций Морса. Пусть
$\D^0\subset\D^\pm$ --- пространство диффео\-мор\-физ\-мов,
гомотопных $\id_M$ в классе гомеоморфизмов.
Снабдим $C^\infty$-топологией пространства $F,\FF^1,\D^0$,
см.~\cite[\S4]{kp1}. Из результатов~\cite {S,EE,EE0} следует, что
имеется гомотопическая эквивалентность
 \begin{equation} \label {eq:EE}
 \D^0 \aapprox R_{\D^0},
 \end{equation}
где $R_{\D^0}$ --- одно из многообразий $SO(3)$ (при $M=S^2$),
$T^2$ (при $M=T^2$) и точка (при $\chi(M)<0$).

Предположим, что количество отмеченных критических точек $\hat p+\hat
q+\hat r >\chi(M)$. Пусть
$$
 \KK:=\KK_{p,q,r;\hat p,\hat q,\hat r} \ \subset \
 \MM:=\MM_{p,q,r;\hat p,\hat q,\hat r}
$$
-- комплекс оснащенных функций Морса и содержащее его $3q$-мерное
многообразие (см.~\cite[\S4]{KP3}). Согласно \cite[\S4]{KP3addmz},
имеется гомеоморфизм $\overline{\Ev}\:\FF^1/\D^0\to\MM$. Рассмотрим
универсальное пространство модулей
$\KK^\infty:=\overline{\Ev}(\FF^0/\D^0)\subset\MM$ специальных
оснащенных функций Морса.

\begin{Thm} \label {thm:KP4addvest}
Пусть выполнены условия теоремы~{\rm\ref {thm:FF0}}, и количество
отмеченных критических точек
 \begin{equation} \label {eq:main}
\hat p+\hat q+\hat r >\chi(M).
 \end{equation}
Имеются гомотопические эквивалентности и гомеоморфизм $F \sim F^1 \sim \FF^0
 \approx \D^0\times\KK^\infty
 \sim R_{\D^0}\times\KK$.
\end{Thm}

\begin{Def} \label{def:equiv}
Функции Морса $f,g\in F$ назовем {\em изотопными} ($f\sim_\isot g$),
если найдутся такие диффеоморфизмы $h_1\in\D^0$ и
$h_2\in\Diff^+(\RR)$, что $f=h_2\circ g\circ h_1$ и $h_1$ сохраняет
нумерацию отмеченных критических точек.
Множество функций из $F^1$, изотопных $f$, обозначим через $[f]$.
\end{Def}

\begin{Lem} \label {lem:FF0}
Имеется $\D^0$-эквивариантный гомеоморфизм
$\FF^0\approx\D^0\times\KK^\infty$.
При этом $\KK^\infty$ является сильным деформационным ретрактом $\MM$, а $\KK$
-- сильным деформационным ретрактом $\KK^\infty$.
\end{Lem}

\Proof Первое утверждение леммы следует из
существования $\D^0$-эквивариантного гомеоморфизма
$\FF^1\approx\D^0\times\MM$, согласованного с $\overline{\Ev}$ (см.\
\cite[\S4]{KP3addmz}). В силу теоремы \ref {thm:FF0},
отображение включения
$\FF^0/\D^0\approx\KK^\infty\hookrightarrow\MM\approx\FF^1/\D^0$
является гомотопической эквивалентностью.

Согласно \cite {KP3}, имеются замкнутое покрытие
$\KK=\cup_{[f]}\DD_{[f]}$ и открытое покрытие
$\MM=\cup_{[f]}\MM_{\succeq[f]}$, такие что $\DD_{[f]}\subset\MM_{\succeq[f]}$,
а также имеются выпуклые множества
$D_{[f]}\subset S_{\succeq[f]}\subset H^0_{[f]}\cong\RR^q$,
$U_{[f]}\subset U_{[f]}^\infty\subset H^1_{[f]}
\cong\RR^{2q}$, такие что отображение включения $\DD_{[f]}\hookrightarrow\MM_{\succeq[f]}$ есть
композиция
 $\DD_{[f]}\approx(D_{[f]}\times U_{[f]})/\tilde\Gamma_{[f]}\subset(S_{\succeq[f]}\times U_{[f]}^\infty)/\tilde\Gamma_{[f]}\approx\MM_{\succeq[f]}$,
где $\tilde\Gamma_{[f]}$ -- группа, действующая свободно, дискретно и
покомпонентно на $S_{\succeq[f]}\times U_{[f]}^\infty$ (более точные
определения см.\ в \cite[\S4]{KP3}).
Из определения ``вычисляющего'' отображения $\Ev\:\FF^1\to\MM$ (см.\
\cite[\S4]{KP3addmz}) следует, что
$\KK^\infty=\cup_{[f]}\DD^\infty_{[f]}$, где
$\DD^\infty_{[f]}\approx(D_{[f]}\times U_{[f]}^\infty)/\tilde\Gamma_{[f]}$.
Отсюда $(\MM,\KK^\infty)$ -- пара полиэдров (а потому корасслоение по
теореме Борсука, см.~\cite [\S5.5]{FF}). Поэтому $\KK^\infty$ --
сильный деформационный ретракт многообразия $\MM$ (см.\ \cite[Гл.\,1,
\S4]{Sp}).

Сильная деформационная ретракция для пары $\KK\subset\KK^\infty$
получается из сильных деформационных ретракций для пар соседних
пространств в цепочке
$\KK=\KK\cup\KK^\infty_{-1}\subset\KK\cup\KK^\infty_0\subset\KK\cup\KK^\infty_1\subset\ldots\subset\KK\cup\KK^\infty_{q-1}=\KK^\infty$,
где $\KK^\infty_{k}:=\cup_{\dim D_{[f]}\le
k}\DD^\infty_{[f]}\subset\KK^\infty$. Сильную деформационную
ретракцию для пары
$\KK\cup\KK^\infty_{k-1}\subset\KK\cup\KK^\infty_k$ определим с
помощью сильных деформационных ретракций для пар $((D_{[f]}\times
U_{[f]}) \cup (\d D_{[f]}\times U_{[f]}^\infty))/\tilde\Gamma_{[f]} \subset
(D_{[f]}\times U_{[f]}^\infty)/\tilde\Gamma_{[f]}$, таких что $\dim
D_{[f]}=k$ (см.\ \cite[\S4]{KP3}).
\qed

\PProof{теоремы \ref {thm:KP4addvest}} По теореме \ref {thm:FF0} и
лемме \ref {lem:FF0} выполнено
$\FF^1\sim\FF^0\approx\D^0\times\KK^\infty\sim\D^0\times\KK$.
С учетом (\ref {eq:EE}) и того, что забывающее отображение и
отображение включения $\FF^1\to F^1\hookrightarrow F$
являются гомотопическими эквивалентностями (см.\ \cite[теорема
2.5]{kp1}), получаем теорему~\ref {thm:KP4addvest}. \qed

\medskip
Автор приносит благодарность С.А.\ Мелихову, Д.А.~Пермякову и
А.Т.~Фоменко за полезные замечания и обсуждения.

Работа частично поддержана грантом РФФИ \No~10--01--00748-а,
грантом про\-граммы ``Ведущие научные школы РФ'' НШ-3224.2010.1,
грантом программы ``Развитие научного потенциала высшей школы'' РНП
2.1.1.3704 <<Современная дифференциальная геометрия, топология и
приложения>> и грантом ФЦП <<Научные и научно-педагогические кадры
инновационной России>> (контракты \No~02.740.11.5213 и \No~14.740.11.0794).


\end{document}